\documentclass{ieeeaccess}
\usepackage{cite}
\usepackage{amsmath,amssymb,amsfonts,array}
\usepackage{ntheorem}
\theoremheaderfont{\it}
\theoremseparator{.}
\usepackage{algorithmic,algorithm,url}
\usepackage{graphicx}
\usepackage{textcomp}
\usepackage[T1]{fontenc}
\usepackage{xcolor,soul}

\def\BibTeX{{\rm B\kern-.05em{\sc i\kern-.025em b}\kern-.08em
    T\kern-.1667em\lower.7ex\hbox{E}\kern-.125emX}}

\newtheorem{thm}{Theorem}
\newtheorem{prop}[thm]{Proposition}
\newtheorem{lem}[thm]{Lemma}
\newtheorem{cor}[thm]{Corollary}
\newtheorem{rmk}[thm]{Remark}
\newtheorem{defn}[thm]{Definition}

\renewcommand{\IEEEQED}{\IEEEQEDopen}

\makeatletter
\newcommand{\SWITCH}[1]{\STATE \textbf{switch} (#1)}
\newcommand{\ENDSWITCH}{\STATE \textbf{end switch}}
\newcommand{\CASE}[1]{\STATE \textbf{case} #1\textbf{:} \begin{ALC@g}}
\newcommand{\ENDCASE}{\end{ALC@g}}

\newcommand{\DEFAULT}{\STATE \textbf{default:} \begin{ALC@g}}
\newcommand{\ENDDEFAULT}{\end{ALC@g}}
\newcommand{\DEFAULTLINE}[1]{\STATE \textbf{default:} }
\makeatother

\begin{document}
\doi{}

\title{Powers of large matrices on GPU platforms to compute the Roman domination number of cylindrical graphs}

\author{\uppercase{J.A. Mart\'inez}\authorrefmark{1,}\authorrefmark{3}, \uppercase{E.M. Garz\'on}\authorrefmark{1,}\authorrefmark{3} and
\uppercase{M.L. Puertas\authorrefmark{2,}\authorrefmark{3}}}
\address[1]{Department of Computer Science, Universidad de Almer\'ia, 04120 Almer\'ia, Spain }
\address[2]{Department of Mathematics, Universidad de Almer\'ia, 04120 Almer\'ia, Spain}
\address[3]{Agrifood Campus of International Excellence (ceiA3), Universidad de Almer\'ia, 04120 Almer\'ia, Spain}
\tfootnote{This work was supported in part by grants RTI2018-095993-B-I00 and PID2019-104129GB-I00/AEI/10.13039/501100011033 of the Spanish Ministry of Science and Innovation.}

\markboth
{J.A. Mart\'inez \headeretal: Powers of large matrices on GPU platforms to compute the Roman domination number of cylindrical graphs}
{J.A. Mart\'inez \headeretal: Powers of large matrices on GPU platforms to compute the Roman domination number of cylindrical graphs}

\corresp{Corresponding author: J.A. Mart\'inez (e-mail: jmartine@ual.es).}

\begin{abstract}

The Roman domination in a graph $G$ is a variant of the classical domination, defined by means of a so-called Roman domination function $f\colon V(G)\to \{0,1,2\}$ such that if $f(v)=0$ then, the vertex $v$ is adjacent to at least one vertex $w$ with $f(w)=2$. The weight $f(G)$ of a Roman dominating function of $G$ is the sum of the weights of all vertices of $G$, that is, $f(G)=\sum_{u\in V(G)}f(u)$. The Roman domination number $\gamma_R(G)$ is the minimum weight of a Roman dominating function of $G$. In this paper we propose algorithms to compute this parameter involving the $(\min,+)$ powers of large matrices with high computational requirements and the GPU (Graphics Processing Unit) allows us to accelerate such operations. {Specific routines have been developed to efficiently compute the $(\min ,+)$ product on GPU architecture, taking advantage of its computational power. These algorithms allow us to }compute the Roman domination number of cylindrical graphs $P_m\Box C_n$ i.e., the Cartesian product of a path and a cycle, in cases $m=7,8,${$9$} $ n\geq 3$ and $m\geq ${$10$}$, n\equiv 0\pmod 5$. Moreover, we provide a lower bound for the remaining cases $m\geq${$10$}$, n\not\equiv 0\pmod 5$.

\end{abstract}


\titlepgskip=-15pt

\maketitle

\section{Introduction}

The efficient location of resources in a network is a well known optimization problem that is usually approached by using graphs. The domination parameters in graphs play a central role in such problems since they can represent a wide variety of additional properties required in the distribution of resources. A dominating set in a graph $G$ is a vertex subset $S$ such that every vertex not in $S$ has at least one neighbor in it. The domination number of $G$ is the cardinal of a minimum dominating set. Applications of this parameter and some of its variations to the optimal location of radio stations or land surveying sensors can be
found in~\cite{Hay1998}. Moreover, linear algorithms using domination parameters for the resource allocation in trees are studied in~\cite{Hedetniemi1986}.

In this paper we focus on the Roman domination that models a classical optimization problem (see~\cite{Stewart1999}). The Roman Emperor Constantine the Great, in the 4th century AD, ordered that no legion be sent out of its usual place if such place was left unprotected. Therefore, a pair of legions must be placed in some locations so that one of them could be sent to an adjacent one. Meanwhile, locations with just one legion just protect themselves. The goal is selecting the appropriate locations to place either one or two legions in order to minimize the needed forces.

Following this approach, the Roman domination in graphs was introduced in~\cite{Cockayne2004}. A {\it Roman dominating function} in a graph $G$ is a function $f\colon V(G)\to \{0,1,2\}$  such that every vertex $v$ with $f(v)=0$ is adjacent to at least a vertex $w$ satisfying $f(w)=2$. The weight of a Roman dominating function is $f(G)=\sum_{u\in V(G)}f(u)$. The minimum weight of a Roman dominating function of $G$ is the {\it Roman domination number} $\gamma_R(G)$. We denote $S^f_i=\{v\in V(G)\colon f(v)=i\}$  (and we will omit $f$ if there is no confusion). Therefore, $f(G)=|S^f_1|+2|S^f_2|$. We say that vertices in $S_1$ just dominate themselves while every vertex $w\in S_2$ dominates itself and its neighborhood $N(w)$.

The Cartesian product of two graphs $G\Box H$ is the graph with vertex set $V(G)\times V(H)$ such that two vertices $(g_1,h_1), (g_2,h_2)$ are adjacent in $G\Box H$ if either $g_1=g_2$ and $h_1, h_2$ are adjacent in $H$, or $g_1,g_2$ are adjacent in $G$ and $h_1=h_2$ (see~\cite{Imrich2000}). The Roman domination number remains unknown for general Cartesian product graphs while the Roman domination number of particular cases of Cartesian product of paths and cycles have been computed. The Roman domination number of the Cartesian product of two paths has been recently obtained in~\cite{Rao2019} meanwhile the problem is still open for the Cartesian product of a path and a cycle and the Cartesian product of two cycles. In both of them, solutions for small cases have been provided by using an algorithmic approach (see~\cite{Pavlic2012}).

Regarding the computational complexity of theses problems, the computation of the domination number is NP-complete in general
graphs (see~\cite{Gar1979}). Moreover, it remains NP-complete when restricted to bipartite or chordal graphs and it is polynomial in a few graphs classes such as trees or interval graphs (see~\cite{Hay1998}).

In a similar way, the computation of the Roman domination number is also an NP-complete problem, even in bipartite, planar or chordal graphs and it can be computed in linear time in
trees (see~\cite{Cockayne2004}). In the same paper, the authors also stated that there is a $2\log n$ approximation algorithm for the Roman domination number. More recently, a $2(1+\log(\Delta+1))$-approximation algorithm, with $\Delta$ the maximum degree of the graph, has been presented in~\cite{Pada2000} to find a Roman domination function with minimum weight.

The great interest to express graph algorithms in terms of tropical algebra operations is well known~\cite{Kepner2011}. From a computational point of view, this approach involves several challenges due to the high dimension of the matrices that take part in such algorithms. This way, many works have focused on taking advantage of the sparsity of the matrices and the regularity of specific graphs to reduce the complexity of the corresponding matrix computations \cite{Dobo1990,Duan2009,Felz2011}; the goals of other works are the optimal computational implementations of the primitive matrix operations related to this field to exploit modern multicore and GPU platforms~\cite{Humayun2016,Kepner2015,Yang2020}.

In this work, the study of the Roman domination in graphs relies on the computation of matrix powers. It is expressed as a $(\min,+)$ product sequence which starts with sparse matrices which are filled as new products are computed. Therefore, the best option is to compute the matrix power with dense data structures. However, it implies a very high computational complexity in terms of run-time and memory requirements since the size of matrices strongly increases with the dimensions of the cylindrical graphs.

The cylinder $P_m\Box C_n$ is the Cartesian product of the path with $m$ vertices $P_m$ and the cycle with $n$ vertices $C_n$, that is, the graph with vertex set $V(P_m)\times V(C_n)$  such that two vertices $(u_1,v_1), (u_2,v_2)$ are adjacent in $P_m\Box C_n$ if $u_1=u_2$ and $v_1, v_2$ are adjacent in $C_n$ or $u_1,u_2$ are adjacent in $P_m$ and $v_1=v_2$. As we have said, there is no general formula for the Roman domination number in this graph family and just exact values of $\gamma_R(P_m\Box C_n)$ with $2\leq m\leq 6$ or $2\leq n\leq 8$ are known (see~\cite{Pavlic2012}).

In this paper we provide an algorithm involving powers of large matrices, that allows us to compute the exact values of $\gamma_{R}(P_7\Box C_n)$ and  $\gamma_{R}(P_8\Box C_n)$. To this end, we have developed a GPU version of $(\min,+)$ powers of dense matrices relying on the routine MatrixMult of~\cite{nvidia2020}. Moreover, graphs of higher dimensions could be studied on GPUs with larger memory capacity than used in the experimental study. Finally, we will use a modification of the algorithm to obtain a lower bound of $\gamma_{R}(P_m\Box C_n)$, for $m\geq 9$, that gives the exact value of this parameter if $n\equiv 0\pmod 5$.

\section{Related work}

The Roman domination in graphs has been widely studied since it was formally defined in 2004 (see~\cite{Cockayne2004}) and several hundreds of papers about it can be found in literature. We quote some references as an example: for instance, general upper bounds of the Roman domination number in terms of the number of vertices of the graphs were first obtained in~\cite{Chambers2009} and this problem has also been studied in~\cite{Jafari2019_1,Liu2012}.

A different point of view is the study of the relationship between the Roman domination number and other domination parameters, for instance with the domination number in~\cite{Favaron2009} or with variations of the Roman domination itself, such as the Roman-$\{2\}$-domination in~\cite{Martinez2019}. Moreover, the Roman domination number can be defined in directed graphs and in~\cite{Hao2019} the relationship between the domination number and the Roman domination number in such graphs is considered. The algorithmic point of view has also been studied, as an example in~\cite{Liedloff2005} authors showed that the Roman domination number of cographs and interval graphs can be computed in linear time.

Furthermore, a number of variations of the original definition have recently been analyzed. As an example, the Italian domination number is computed in some Cartesian product graphs with algorithmic procedures in~\cite{Gao2019}; the sum and the product of the double Roman domination numbers of a graph and its complement are considered in~\cite{Jafari2019_2}; the Roman-$\{2\}$-bondage number is computed in~\cite{Moradi2020} for some graph families such as paths, cycles, complete bipartite graphs, trees, uniclyclic graphs and planar graphs; the total double Roman domination number is studied in~\cite{Sao2019}, where some general upper bounds in terms of the number of vertices and the maximum degree are obtained.

It is well-known that the domination properties are difficult to handle in the Cartesian product graphs and computing such parameters in them has attracted attention since the Vizing Conjecture, which remains open, was formulated in 1968 (see~\cite{Vizing1968}). As an example, the problem of computing the domination number of the Cartesian product of two paths was open for almost thirty years and it was finally solved in~\cite{Goncalves2011}, where authors provide the value of this parameter in terms of the number of vertices of both paths. In previous works, upper and lower bounds can be found, as in~\cite{Cockayne1985,Gravier1997,Jacobson1983}. Moreover, in~\cite{Jacobson1983,Spalding1998} the exact values of some small cases were computed. This problem is still open in other Cartesian products of two graphs involving paths, cycles or more general graphs. Also, some partial results about the domination number of the Cartesian product of $k$ paths can be found in~\cite{Gravier1997} and for the case of $k$ cycles in~\cite{Klavzar1995}.

Among the techniques used to compute domination parameters in Cartesian product graphs, an algorithmic approach using matrix powers has provided significant results in cases involving paths and cycles. The final paper for the domination number of the Cartesian product of two paths~\cite{Goncalves2011} uses this technique, which was presented in~\cite{Klavzar1996} for fasciagraphs and rotagraphs, of which Cartesian products of paths and cycles are particular cases. Later, in~\cite{Pavlic2013} the authors followed these ideas to compute the domination number of the Cartesian product of two cycles and the Cartesian product of a path and cycle, but just in some small cases. The same authors have adapted the technique to obtain the Roman domination number of some Cartesian product graphs involving small paths and/or cycles in~\cite{Pavlic2012}. Recently, in~\cite{Rao2019} the Roman domination number of the Cartesian product of two paths has been completely computed following the same ideas.

\section{Preliminary results}

In this section we present the needed tools to provide the algorithms that we will use to compute the Roman domination number in selected cylinders.

The first tool is the $(\min,+)$ matrix algebra over the semi-ring $\mathcal{P}=(\mathbb{R}\cup\{\infty\}, \min, +, \infty, 0)$ of tropical numbers in the minimum convention (see~\cite{Pin1998}). The $(\min, +)$ matrix multiplication $\bigotimes$ is defined by $C=A \bigotimes B$, being the matrix where for all $i,j$, $c_{i,j}=\min \limits_{k}(a_{i,k}+b_{k,j})$.

Moreover, the $(\min,+)$ product of a matrix $A$ and $\alpha\in \mathbb{R}\cup\{\infty\}$ is defined by $(\alpha\bigotimes A)_{i,j}=\alpha+a_{i,j}$.
Therefore, $(A\bigotimes (\alpha \bigotimes B))_{i,j}=\min \limits_{k}(a_{i,k}+(\alpha \bigotimes B)_{k,j})=\min \limits_{k}(a_{i,k}+(\alpha+b_{k,j}))=
\min \limits_{k}(\alpha+(a_{i,k}+b_{k,j}))=\alpha+\min \limits_{k}(a_{i,k}+b_{k,j})=\alpha + (A\bigotimes B)_{i,j}$. Hence $A\bigotimes (\alpha\bigotimes B)=\alpha \bigotimes(A\bigotimes B)$.

The second ingredient that we need is the following result (see~\cite{Carre1979}) that we quote from~\cite{Klavzar1996}. We just need the particular case related to the tropical semi-ring.

Let $G$ be a digraph with $V(G)= \{v_1, v_2, \dots, v_s\}$ together with a labeling function $\ell$ which assigns an element of $\mathcal{P}$ to every arc of $G$. A path of length $k$ in $G$ is a sequence of $k$ consecutive arcs $Q=(v_{i_0}v_{i_1})(v_{i_1}v_{i_2})\dots (v_{i_{k-1}}v_{i_k})$ and $Q$ is \emph{a closed path} if $v_{i_0}=v_{i_{k}}$. The labeling $\ell$ can be easily extended to paths
$$\ell(Q) = \ell(v_{i_0}v_{i_1})+\ell(v_{i_1}v_{i_2})+\dots +\ell(v_{i_{k-1}}v_{i_{k}}).$$
\begin{thm}\label{thm:carre}
Let $S_{ij}^k$ be the set of all paths of length $k$ from $v_i$ to $v_j$ in $G$ and let $A(G)$ be the
matrix defined by
$$A(G)_{ij} =
\left\{
\begin{array}{ll}
\ell(v_i,v_j) & \text{if }(v_i,v_j) \text{ is an arc of } G,\\
\infty & \text{otherwise.}
\end{array}
\right.
$$
If $A(G)^k$ is the $k$-th $(\min, +)$ power of $A(G)$, then
$$(A(G)^k)_{ij}=\min \{\ell(Q)\colon Q\in S_{ij}^k\}.$$
\end{thm}

Finally, we will also use the following lemma about the {$(\min,+)$} matrix multiplication. It is a standard argument and we include its proof here for the sake of completeness.

\begin{lem}\label{lem:recurrence}
Let $A$ be a square matrix and suppose that there exist natural numbers $n_0,a,b$ such that $A^{n_0+a}=b\bigotimes A^{n_0}$ then, $A^{n+a}=b\bigotimes A^{n}$, for every $n\geq n_0$.
\end{lem}

\begin{IEEEproof}
We proceed by induction. By hypothesis, $A^{n_0+a}=b\bigotimes A^{n_0}$. Let $n\geq n_0$ be such that $A^{n+a}=b\bigotimes A^{n}$ then, $A^{(n+1)+a}=A\bigotimes A^{n+a}=A\bigotimes(b\bigotimes A^{n})=b\bigotimes(A\bigotimes A^{n})=b\bigotimes A^{n+1}$, as desired.
\end{IEEEproof}

\section{Computation of $\gamma_{R}(P_7\Box C_{\lowercase{n}})$ and $\gamma_{R}(P_8\Box C_{\lowercase{n}})$}\label{sec:small}

In this section we provide an algorithm to compute the exact values of $\gamma_{R}(P_m\Box C_{\lowercase{n}})$, for $m=7,8$, by using Theorem~\ref{thm:carre}. Our approach follows the ideas in~\cite{Pavlic2012,Rao2019}. In the first paper, authors compute the Roman domination number of $P_m\Box C_n$, with $m\leq 6$, and we use a modification of their algorithm, following the techniques shown in the second reference to compute the following two cases.

Our algorithm constructs a matrix for each value of $m$ and computes some $(\min, +)$ powers of it. The sizes of such matrices grow exponentially with $m$ and it is not expected that it could compute cases much larger than $m=8$, even by using additional computing resources. However, these small cases play an important role in formulating a conjecture about the behavior of the general case. Although it does not happen for very small $m$, a regular behavior could appear as in the case of the Cartesian product of two paths. This is {the reason for} our interest to compute as many small cases as possible.

\subsection{Theoretical results}\label{subsec:th}

First of all, we encode the vertex set of $P_m\Box C_n$. We say that it has $m$ rows and $n$ columns, each row being a cycle with $n$ vertices and each column being a path with $m$ vertices (see Figure~\ref{fig:cylinder}). The $m$ rows are numerated from top to bottom and the $n$ columns from left to right. We will consider that the last column is the previous column of the first one (and the first column follows the last one).

\Figure[!h]()[width=0.28\textwidth]{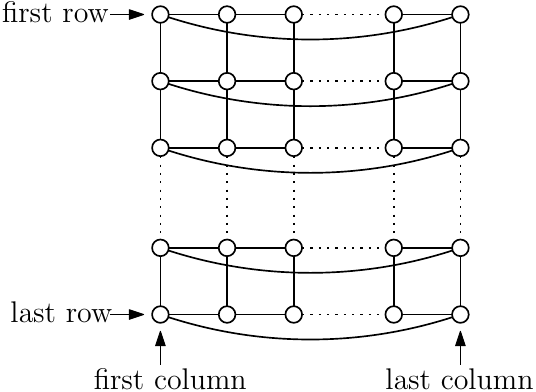}
{The cylinder $P_m \Box C_n$. \label{fig:cylinder}}

Let $f$ be a Roman dominating function of $P_m \Box C_n$, with associated vertex partition $S_i\! =\!\{v\!\in \! V(P_m \Box C_n)\! \colon \!f(v)\!=\!i\}$, for $i=0,1,2$. We assign a label to each vertex in the following way.

\begin{itemize}
\item $v=a$ if $v\in S_2$,
\item $v=b$ if $v\in S_1$,
\item $v=c$ if $v\in S_0$ and it has at least one neighbor in $S_2$ in its column or in the previous one,
\item $v=d$ if $v\in S_0$ and it has no neighbor in $S_2$ in its column nor in the previous one.
\end{itemize}

Now each column is a word of length $m$ in the alphabet $\{a,b,c,d\}$ and the function $f$ can be uniquely identified with a sequence of $n$ consecutive words. Note that not any word can appear associated to such functions because the definition of the labeling implies that letter sequences $ad, da$ are not possible.

\begin{rmk}\label{rem:reduced_set}
If a word associated to a Roman dominating function contains any of the sequences $ab,ba,bb$, we can replace them with $ac,ca,ac$ respectively to obtain another Roman dominating function with, at most, the same weight.
\end{rmk}

Bearing in mind these considerations, we pose the following definition.

\begin{defn}
A correct word of length $m$ in the alphabet $\{a,b,c,d\}$ is a sequence of $m$ letters not containing $ad,da,ab,ba,bb$.
\end{defn}

Note that some correct words cannot appear together associated to a Roman dominating function, again because of the definition of the labeling. We now list the needed conditions for a pair of correct words $\mathbf{p}=(p_1,\dots p_m)$ and $\mathbf{q}=(q_1, \dots q_m)$ (we will denote words in bold font and their letter in italic font), in order to make it possible that {$\mathbf{p}$} can follow {$\mathbf{q}$}. Note that we also avoid combinations $ab,ba,bb$ in the same row (see Remark~\ref{rem:reduced_set}), in order to reduce the number of suitable sequences of words.

\begin{enumerate}
\item conditions for the first row

\begin{itemize}
\item if $q_1=a$, then $p_1=a$ or $p_1=c$
\item if $q_1=b$, then $\big( p_1=c \text{{ and} } p_2=a \big )$ or $p_1=d$
\item if $q_1=c$, then $p_1=a$ or $p_1=b$ or $\big( p_1=c \text{{ and} }  p_2=a \big)$ or $p_1=d$
\item if $q_1=d$, then $p_1=a$
\end{itemize}

\item conditions for the intermediate rows $2\leq i\leq m-1$

\begin{itemize}
\item if $q_i=a$, then {$p_i$}$=a$ or {$p_i$}$=c$
\item if $q_i=b$, then $\big( p_i=c \text{{ and} } p_{i-1}=a \big)$ or $\big( p_i=c \text{{ and} } p_{i+1}=a\big)$ or $p_1=d$
\item if $q_i=c$, then $p_i=a$ or $p_i=b$ or $\big( p_i=c \text{{ and} } p_{i-1}=a\big)$ or $\big( p_i=c \text{{ and} } p_{i+1}=a\big)$ or $p_1=d$
\item if $q_i=d$, then $p_i=a$
\end{itemize}

\item conditions for the last row
\begin{itemize}
\item if $q_m=a$, then $p_m=a$ or $p_m=c$
\item if $q_m=b$, then $\big( p_m=c \text{{ and} } p_{m-1}=a\big)$ or $p_m=d$
\item if $q_m=c$, then $p_m=a$ or $p_m=b$ or $\big( p_m=c \text{{ and} } p_{m-1}=a\big)$ or $p_m=d$
\item if $q_m=d$, then $p_m=a$
\end{itemize}

\end{enumerate}

\begin{rmk} \label{rem:min}
Denote by $\mathcal{S}$ the set of all the sequences {$\mathbf{p_1},\mathbf{p_2}\dots \mathbf{p_n}$} of $n$ correct words of length $m$ such that $\mathbf{p_1}$ can follow $\mathbf{p_n}$ and $\mathbf{p_{i+1}}$ can follow $\mathbf{p_i}$, for $1\leq i\leq n-1$. From Remark~\ref{rem:reduced_set}, the set $\mathcal{R_S}$ of Roman dominating functions associated to the sequences in $\mathcal{S}$ contains at least one function with minimum weight.
\end{rmk}

In order to apply Theorem~\ref{thm:carre}, we consider the digraph $G$ whose vertex set is $V(G)=\{ \text{correct words of length } m $ $\text{in the alphabet } \{a,b,c,d\} \}$, and such that there is an arc from the word $\mathbf{q}$ to the word $\mathbf{p}$ if $\mathbf{p}$ can follow $\mathbf{q}$. Moreover, we define the following labeling function $\ell$ which assigns to every arc of $G$ an element of the semi-ring of tropical numbers $\mathcal P$: $\ell(\mathbf{q},\mathbf{p})=2 \mathbf{p}(a)+\mathbf{p}(b)$, where $\mathbf{p}(a)=$number of $a's$ in the word $\mathbf{p}$ and $\mathbf{p}(b)=$number of $b's$ in the word $\mathbf{p}$.

\begin{prop}\label{prop:weight}
Let $f\in \mathcal{R_S}$ be a Roman dominating function of $P_m\Box C_n$ and let $\mathbf{p_1} \mathbf{p_2} \dots \mathbf{p_n}$ be the sequence of correct words associated to $f$. Then,
$Q=(\mathbf{p_1} \mathbf{p_2})(\mathbf{p_2} \mathbf{p_3})\dots (\mathbf{p_{n-1}} \mathbf{p_n})(\mathbf{p_n} \mathbf{p_1})$ is a closed path in the digraph $G$ such that $\ell (Q)=f(P_m\Box C_n)$, the weight of $f$.
\end{prop}

\begin{IEEEproof}
Clearly $Q$ is a closed path in the digraph $G$, because there is an arc from $\mathbf{p_i}$ to $\mathbf{p_{i+1}}$ (for $1\leq i\leq n-1$) and from $\mathbf{p_n}$ to $\mathbf{p_1}$. Moreover, $\ell (Q)= \ell(\mathbf{p_1} \mathbf{p_2})+\ell(\mathbf{p_2} \mathbf{p_3})+\dots +\ell(\mathbf{p_{n-1}} \mathbf{p_n})+\ell(\mathbf{p_n} \mathbf{p_1})=(2 \mathbf{p_2}(a)+\mathbf{p_2}(b))+(2 \mathbf{p_3}(a)+\mathbf{p_3}(b))+\dots +(2 \mathbf{p_n}(a)+\mathbf{p_n}(b))+(2 \mathbf{p_1}(a)+\mathbf{p_1}(b))=2(\mathbf{p_1}(a)+\mathbf{p_2}(a)+\dots \mathbf{p_n}(a))+(\mathbf{p_1}(b)+\mathbf{p_2}(b)+\dots \mathbf{p_n}(b))=2|S^f_2|+|S^f_1|=f(P_m\Box C_n)$.
\end{IEEEproof}

\begin{cor}\label{cor:diagonal}
Let $A(G)$ be the matrix defined by
$$A(G)_{\mathbf{qp}} =
\left\{
\begin{array}{ll}
\ell(\mathbf{q},\mathbf{p}) & \text{if }(\mathbf{qp}) \text{ is an arc of } G,\\
\infty & \text{otherwise.}
\end{array}
\right.
$$
Then $\min _\mathbf{p} (A(G)^n)_{\mathbf{pp}}=\gamma _R(P_m\Box C_n)$.
\end{cor}

\begin{IEEEproof}
From Theorem~\ref{thm:carre} we know that $(A(G)^n)_{\mathbf{pp}}=\min \{\ell(Q)\colon Q\in S_{\mathbf{pp}}^n\}$, where $S_{\mathbf{pp}}^n$ is the set of all closed paths of length $n$ from $\mathbf{p}$ to $\mathbf{p}$ in $G$. Clearly, the closed path
$Q=(\mathbf{p} \mathbf{p_2})\dots (\mathbf{p_{n-1} p_{n}})(\mathbf{p_{n} p})$ belongs to $S_{\mathbf{pp}}^n$ if and only if the sequence of words
$\mathbf{p,  p_2,\dots p_n}$ belong to the set $\mathcal{S}$. Therefore, from Remark~\ref{rem:min} and Proposition~\ref{prop:weight},
$\min_\mathbf{p}(A(G)^n)_{\mathbf{pp}}=\min_\mathbf{p} (\min\{\ell(Q)\colon Q\in S_{\mathbf{pp}}^n\})=\min_\mathbf{p}(\min \{f(P_m\Box C_n)\colon f \in \mathcal{R_S}, \mathbf{p} \text{ is the first column} \})=\gamma _R(P_m\Box C_n)$.
\end{IEEEproof}
\par\bigskip

The Corollary above allows us to compute the Roman domination number of a cylinder $P_m\Box C_n$, where both $m$ and $n$ are fixed.
The recurrence argument for $(\min, +)$ matrix powers shown in Lemma~\ref{lem:recurrence} gives the opportunity to compute the Roman domination number of $P_m\Box C_n$, just fixing one of the sizes but not both of them, as we show in the following result.

\begin{prop}\label{pro:difference}
Let $m$ be a natural number and consider the digraph $G$ and the matrix $A(G)$ constructed above. If $A(G)^{n_0+\alpha}=\beta \bigotimes A(G)^{n_0}$ for natural numbers $n_0,\alpha, \beta$, then $\gamma _R(P_m\Box C_{n+\alpha})-\gamma_R(P_m\Box C_n)=\beta$, for $n\geq n_0$.
\end{prop}

\begin{IEEEproof}
If $A(G)^{n_0+\alpha}=\beta \bigotimes A(G)^{n_0}$ then, by using Lemma~\ref{lem:recurrence} we obtain that $A(G)^{n+\alpha}=\beta \bigotimes A(G)^{n}$, for $n\geq n_0$. Moreover, {Corollary}~\ref{cor:diagonal} gives
$\gamma _R(P_m\Box C_{n+\alpha})=\min _\mathbf{p} (A(G)^{n+\alpha})_{\mathbf{pp}}=\min _\mathbf{p} (\beta \bigotimes A(G)^{n})_\mathbf{pp}=\beta +\min _\mathbf{p} (A(G)^{n})_\mathbf{pp}=\beta +\gamma _R(P_m\Box C_{n})$.
\end{IEEEproof}

\par\bigskip

The unique solution of the finite difference equation $\gamma _R(P_m\Box C_{n+\alpha})-\gamma_R(P_m\Box C_n)=\beta, n\geq n_0$, with the boundary values $\gamma_R(P_m\Box C_k), n_0\leq k\leq n_0+(\alpha-1)$, provides $\gamma_R(P_m\Box C_n)$, for $m$ fixed and $n\geq n_0$.

\subsection{The algorithms}\label{subsec:alg}

We can now describe the basic algorithm to compute the Roman domination number of the cylinder $P_m\Box C_n$, {with $m$ and $n$ fixed, by using Corollary}~\ref{cor:diagonal}.

\begin{algorithm}[H]
\caption{Computation of $\gamma_R(P_m\Box C_n)$, for $m,n$ fixed}
\label{alg:small}
\begin{algorithmic}[1]
\REQUIRE {$m,n$ natural numbers}
\ENSURE {$\gamma _R(P_m\Box C_n)$}
\STATE compute all correct words of length $m$
\STATE compute matrix $A(G)$
\STATE compute the $(\min,+)$ matrix power $A(G)^n$
\RETURN $\min _\mathbf{p} (A(G)^n)_{\mathbf{pp}}$
\end{algorithmic}
\end{algorithm}

This algorithm allows us to compute the Roman domination number for cylinders $P_m\Box C_n$ where $m,n$ are small. On the one hand, the size of matrix $A(G)$ is the number of correct words, which is in the order of $4^m$. This size grows exponentially thus, the algorithm is useful only if $m$ is small enough. On the other hand, it is needed to compute the $n-th$ power of the matrix $A(G)$ and, in spite of being a sparse matrix, it becomes dense after a small number of multiplications. This gives that parameter $n$ can not be large so then it is possible to run the algorithm in a short time.

We now present the algorithm to compute the Roman domination number of $P_m\Box C_n$, where just $m$ is fixed. We use Proposition~\ref{pro:difference}, so the algorithm compute the integers $n_0,\alpha, \beta$ that we need to solve the finite difference equation that gives the formula for  $\gamma_R(P_m\Box C_n)$.

\par\medskip

\begin{algorithm}[H]
\caption{Computation of $\gamma_R(P_m\Box C_n)$, for $m$ fixed}
\label{alg:recurrence}
\begin{algorithmic}[1]
\REQUIRE {$m$ a natural number}
\ENSURE {the finite difference equation\\ $\gamma _R(P_m\Box C_{n+\alpha})-\gamma_R(P_m\Box C_n)=\beta$, for $n\geq n_0,$\\
or finite difference equation not found }
\STATE compute all correct words of length $m$
\STATE compute matrix $A(G)$
\STATE compute every $(\min,+)$ power $A(G)^k$, for $k\leq K$ big enough
 \IF {$A^{n_0+\alpha}=\beta \bigotimes A^{n_0}$ for $n_0,\alpha, \beta $ natural numbers}
 \RETURN $n_0, \alpha , \beta$
\ELSE
 \RETURN recurrence not found
 \ENDIF
 \end{algorithmic}
\end{algorithm}

There are some sufficient conditions to ensure that Step 4 in Algorithm~\ref{alg:recurrence} is true (see~\cite{Spalding1998}). However, these theoretical results provide a huge value for $n_0$, in the order of the square of the matrix size, and it is not practical. We have looked for the desired relationship just by checking the powers computed in Step 3, with $K=50$.

\subsection{Computational results}\label{subsec:GPU}

In this subsection we present the details of the implementation in C programming language of the algorithms and the results we have obtained. As we said before, for each $m$ the number $C_m$ of correct words is of the order of $4^m$, which is the number of all words of length $m$ that can be done with $4$ letters.
The first step in both algorithms is the computation of correct words, having a computational complexity of $O(4^m)$. Due to the small values of $m$ considered, this computing time is not relevant.
The second step is the computation of the matrix $A(G)$, whose size $C_m$ grows exponentially with $m$ so, this step is $O((C_m)^2)$. Moreover, this computation has plenty of control flow instructions therefore, it can take better advantage of CPU architecture than GPU architecture.

In Table~\ref{tab:matrix_size} we quote the matrix sizes for $2\leq m\leq 11$, the memory requirements and the computing time of such matrices if $m\leq 9$. {In cases $m=10$ and $11$, we have estimated the memory requirements to store each matrix bearing in mind the matrix sizes and that each entry is an int data (4 bytes).} An Intel(R) Xeon(R) CPU E5-2650 0 @ 2.00GHz with 8 cores (16 threads) and 64 GB of RAM has been used to compute the matrices in both algorithms.

\begin{table}[h!]
\setlength{\extrarowheight}{0.05cm}
\caption{Matrix sizes, memory requirements  and GPU running time of matrix computation line 2 in Algorithm~\ref{alg:small} and Algorithm~\ref{alg:recurrence}}
\setlength{\tabcolsep}{3pt}
\begin{tabular}{|>{\raggedleft\arraybackslash}p{15pt}|>{\raggedleft\arraybackslash}p{55pt}|>{\raggedleft\arraybackslash}p{65pt}|>
{\raggedleft\arraybackslash}p{60pt}|}
\hline
m& Matrix size $C_m$ & Memory size (approx.) & Computing time of matrix $A(G)$\\
\hline
2& 11 & $<1 $MB & $<1 $ {s}\\
\hline
3& 33 &  $<1 $MB & $<1 $ {s}\\
\hline
4& 97 & $<1 $MB& $<1 $ {s}\\
\hline
5& 287 & $<1 $ MB & $<1 $ {s}\\
\hline
6& 848 &  2.7 MB & $<1 $ {s}\\
\hline
7& 2507 & 23.9 MB & {3.21 s}\\
\hline
8& 7411 & {209.5} MB & {30.04 s}\\
\hline
9 & {21909} & 1.78 GB  & {4 min 47 s}\\
\hline
10 & 64769 & estimated 15 GB  & --\\
\hline
11& 191476 & estimated 136 GB& --\\
\hline
\end{tabular}
\label{tab:matrix_size}
\end{table}

Regarding the third step in both algorithms, the complexity of the matrix multiplication operation is $O((C_m)^3)$ so, it consumes most of the running time of both algorithms and it is computed $K$ times to get $A(G)^K$. {The $(\min,+)$ matrix product consists of a modification of the usual product that keeps the complexity.}

To explore the computational requirements of these operations, we have carried out some examples of matrix multiplication by computing one product with an OpenMP CPU implementation (same CPU as before) by using 16 number of threads. In the largest case $m=9$, the matrix $A(G)$ has a size of 21909 and the computation of just one $(\min,+)$ product takes 7 hours and 30 minutes.

Due to the long running time, we have used a GPU to accelerate the matrix powers computation in both algorithms.
Such powers have been carried out by a modification of the
routine MatrixMul, available in the NVIDIA CUDA TOOLKIT 11~\cite{nvidia2020} and described in the CUDA C Programming Guide (see~\cite{nvidia2021}, Chapter 3), to adapt it to
the $(\min,+)$ multiplication. We call this new routine CuMatrixTrop. Such code takes advantage of shared memory to optimize the memory access. We have adapted the kernel code by changing the usual operations in the matrix multiplication: the mapping operation, which is the multiplication, is replaced by the addition and the reduction operation, which is the addition, is replaced by the minimum.

We have used a CUDA block size $32\times 32$ and a padding matrix approach to optimize the performance of CuMatrixTrop. The padding consists of filling the matrix with additional rows and columns of the neutral element, i.e. infinite, until the matrix size is a multiple of 32. These strategies give an efficient exploitation of the GPU, which we have analyzed with the NVIDIA profiling tool {\it nvprof}. In all cases, the achieved occupancy, the global memory load/store efficiency and the warp execution efficiency are near to 100$\%$. The shared memory efficiency is 66$\%$ due to the lower parallelism of the reduction operation.

As we said before, we have run Algorithm~\ref{alg:recurrence} with $K=50$ and this is enough in cases $2\leq m\leq 9$, for which the algorithm returns the integers $n_0, \alpha$ and $\beta$. We have used an NVIDIA Tesla K80 GPU, with 12GB of memory, 13 multiprocessors with 192 cores in each multiprocessor (2496 cores CUDA) and shared memory of 49152 bytes, to compute $A(G)^{k}$ with $k\leq 50$. In Table~\ref{tab:small} we show the results for $2\leq m\leq 9$.

\begin{table}[h!]
\setlength{\extrarowheight}{0.05cm}
\caption{Running times on K80 GPU to compute the matrix powers in Algorithm~\ref{alg:recurrence}}
\setlength{\tabcolsep}{3pt}
\begin{tabular}{|>{\raggedleft\arraybackslash}p{15pt}|>{\raggedleft\arraybackslash}p{60pt}|>{\raggedleft\arraybackslash}p{15pt}|>
{\raggedleft\arraybackslash}p{15pt}|>{\raggedleft\arraybackslash}p{15pt}|}
\hline
m&Computing time of 50 powers& $n_0$ & $\alpha $ & $\beta $ \\
\hline
2& $<1$ {s} & 6& 4 &4\\
\hline
3& $<1$ {s} & 11& 4 &6\\
\hline
4& $<1$ {s} &16& 5& 10\\
\hline
5& {1.19 s} &16& 5& 12\\
\hline
6& {8.42 s} &19 & 5 & 14\\
\hline
7& 1 min 11 {s} & 21 & 5 & 16\\
\hline
8& 11 min 14 {s} & 21 & 5 & 18\\
\hline
9& 76 min 27 s & 22 & 5 & 20\\
\hline
\end{tabular}
\label{tab:small}
\end{table}

Note that it is necessary to save every matrix $A(G)^k$ with $2\leq k\leq 50$, in order to complete steps 4 to 8 in Algorithm~\ref{alg:recurrence}. We have evaluated the profiling of the $A(G)^{50}$ computation based on CuMatrixTrop in terms of GPU running time, disk storage time and CPU-GPU communications. The results are shown in Table~\ref{tab:profiling}.

\begin{table}[h!]
\setlength{\extrarowheight}{0.05cm}
\caption{Profiling results of computing $A(G)^{50}$ in Algorithm~\ref{alg:recurrence}}
\setlength{\tabcolsep}{3pt}
\begin{tabular}{|>{\raggedleft\arraybackslash}p{10pt}|>{\raggedleft\arraybackslash}p{70pt}|>
{\raggedleft\arraybackslash}p{70pt}|>{\raggedleft\arraybackslash}p{65pt}|}
\hline
m&  \ \ \ \ \ \ \ \ CuMatrixTrop ($\%$ over total time) &CPU-GPU communication times ($\%$ over total time) & \ \ \ \ \ \ Disk storage ($\%$ over total time)  \\
\hline
6&  0.30 s (3.53$\%$)&0.02 s (0.24$\%$)& 8.12 s (96.23$\%$) \\
\hline
7& 4.92 s (6.94$\%$)&0.12 s (0.17$\%$) & 66.08 s (92.89$\%$)\\
\hline
8& 139.34 s (20.67$\%$)&1.01 s (0.15$\%$) & 534.66 s (79.18$\%$) \\
\hline
9&  3100.00 s (67.58$\%$)&8.78 s (0.19$\%$) &1487.00 s  (32.23$\%$)\\
\hline
\end{tabular}
\label{tab:profiling}
\end{table}

The disk storage needed in the algorithm consumes a relevant amount of resources. The disk storage is carried out on a Network File System (NFS)
that slows down with the heavy network traffic due to the communication of large matrices. So, this is a bottleneck in our experimental environment. However, the penalties due to the CPU-GPU communications are not significant so, the overload due to GPU exploitation has a negligible impact on performance.

Finally, we have solved the finite difference equations provided by Proposition~\ref{pro:difference}, for $2\leq m\leq ${$9$}, with the parameters shown in Table~\ref{tab:small} and boundary values (that is, the cases $n_0\leq n\leq n_0+(\alpha-1)$) computed with Algorithm~\ref{alg:small}.

The remaining values (cases $3\leq n<n_0$), have also been computed with Algorithm~\ref{alg:small}. We just show the final formula for $m=7,8,${$9$} which are the new cases, for $m\leq 6$ we have obtained the same results as~{\mbox{\cite{Pavlic2012}}}.

\par\bigskip

$\gamma_R(P_7\Box C_n)=
\left\{
\renewcommand*{\arraystretch}{1.5}
\begin{array}{ll}
\big\lceil \frac{16n}{5}\big\rceil & \text{if } n\equiv 0\pmod 5\\
\big\lceil \frac{16n}{5}\big\rceil+1 & \text{otherwise}
\end{array}
\right.
$
\par\bigskip

$\gamma_R(P_8\Box C_n)=
\left\{
\renewcommand*{\arraystretch}{1.5}
\begin{array}{ll}

\big\lceil \frac{18n}{5}\big\rceil & \text{if } n\equiv 0\pmod 5\\

\big\lceil \frac{18n}{5}\big\rceil+1 & \text{if } n\equiv 2,3,4\pmod 5,\\
&  \text{or } n=6\\

\big\lceil \frac{18n}{5}\big\rceil+2 & \text{otherwise}
\end{array}
\right.
$

\par\bigskip

$\gamma_R(P_9\Box C_n)=
\left\{
\renewcommand*{\arraystretch}{1.5}
\begin{array}{ll}
4n & \text{if } n\equiv 0\pmod 5\\
4n+2 & \text{otherwise}
\end{array}
\right.
$

\par\bigskip

\section{A lower bound of $\gamma_{R}(P_{\lowercase{m}}\Box C_{\lowercase{n}})$ and the computation of $\gamma_{R}(P_{\lowercase{m}}\Box C_{\lowercase{5k}})$}\label{sec:bound}

As we have said before, the size of the matrix is a decisive factor when running Algorithm~\ref{alg:small} and Algorithm~\ref{alg:recurrence}, and matrices for values {$m\geq 10$} are too big to allocate {two of them} in the GPU memory, {to perform the matrix product operation}. So we have developed a different approach that we show in this section. We have adapted both algorithms to compute the loss of a Roman dominating function in a cylinder. The loss, originally called wasted domination in~\cite{Guichar2004}, has been used to study several domination parameters in the Cartesian product of two paths (see~\cite{Goncalves2011,Rao2019}). We follow these ideas to compute a lower bound of the Roman domination number for cylinders of large size and with a small computational cost. Although this bound is not tight in general, it provides the exact value when $n\equiv 0\pmod 5$.

We would also like to point out that both versions of the algorithms are complementary. The computation of small cases is needed to conjecture the lower bound, that will later be confirmed by the modified algorithm. Moreover, this second version covers, at least with a lower bound, those cases that are not within the scope of the original algorithm.

Our conjecture for the lower bound of the Roman domination number in cylinder $P_m\Box C_n$ comes from the regular pattern followed by values of $\alpha, \beta$ in Table~\ref{tab:small} for $4\leq m\leq 9$, that is $\alpha =5$ and $\beta =2(m+1)$. Assuming that such pattern remains the same for larger values of $m$, we conjecture that $\big\lceil \frac{2(m+1)n}{5}\big\rceil \leq \gamma_R(P_m\Box C_n)$.

\begin{defn}
Let $f$ be a Roman dominating function of $P_m\Box C_n, m\geq 10$, with weight $f(P_m\Box C_n)$. We define the loss of $f$ as
$L(f)=\frac{5}{2}(f(P_m\Box C_n))-mn$ and we denote the minimum loss of a Roman dominating function in $P_m\Box C_n$ by $\mathcal{L}{(m,n)}$.
\end{defn}

The following lemma is the key result to obtain our lower bound and its proof follows the techniques shown in Section~\ref{sec:small}, with the needed modifications. This proof is quite long and we present it in Appendix~\ref{appendix}, in order to make this section easier to follow.

\begin{lem}\label{lem:bound}
$\mathcal{L}{(m,n)}\geq n$, for $m,n\geq 10$.
\end{lem}

We can now present the announced result with a general lower bound of the Roman domination number in cylinders.

\begin{thm}\label{thm:bound}$\big\lceil \frac{2(m+1)n}{5}\big\rceil \leq \gamma_R(P_m\Box C_n)$, for $m,n\geq 10$.
\end{thm}
\begin{IEEEproof}
$\mathcal{L}(m,n)=\min_f L(f)=\min_f (\frac{5}{2}(f(P_m\Box C_n))-mn)=\frac{5}{2}\gamma_R(P_m\Box C_n)-mn$. This gives that $\gamma_R(P_m\Box C_n)=\frac{2}{5}(mn+\mathcal{L}(m,n))$. Finally, from Lemma~\ref{lem:bound}, $\gamma_R(P_m\Box C_n)=\frac{2}{5}(mn+\mathcal{L}(m,n))\geq \frac{2}{5}(mn+n)=\frac{2(m+1)n}{5}$. Therefore $\big\lceil \frac{2(m+1)n}{5}\big\rceil \leq \gamma_R(P_m\Box C_n)$, as desired.
\end{IEEEproof}
\par\medskip

This bound provides the exact value if $n\equiv 0\pmod 5$.

\begin{cor}
If $m,n\geq ${$10$} and $n\equiv 0\pmod 5$, then $\gamma_R(P_m\Box C_{n}) =\frac{2(m+1)n}{5}\cdot $
\end{cor}

\begin{IEEEproof}
In Figure~\ref{fig:five} we show a Roman dominating function $f$ in {$P_{10}$}$\Box C_{10}$: black vertices have image $2$, grey ones have image $1$ and white ones have image $0$. This regular construction can be repeated in any cylinder with $m,n\geq ${$10$} and $n\equiv 0\pmod 5$ and it has $\frac{n}{5}$ black vertices in each row and $\frac{n}{5}$ grey vertices in the first and the last rows. Therefore, its weight is $f(P_m\Box C_{n})=2m \frac{n}{5}+2\frac{n}{5}=\frac{2(m+1)n}{5}$. This means that $\gamma_R(P_m\Box C_{n})\leq \frac{2(m+1)n}{5}$ and Theorem~\ref{thm:bound} gives the desired equality.
\end{IEEEproof}

\Figure[!h]()[width=0.3\textwidth]{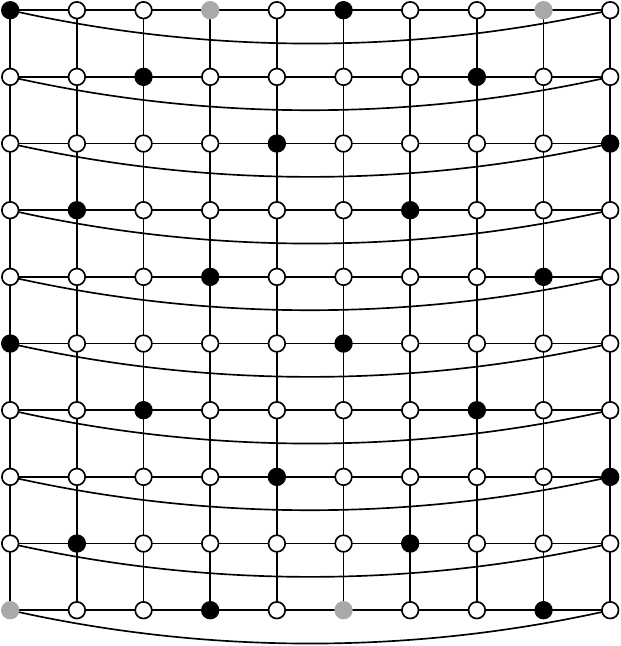}
{Roman domination in {$P_{10}$}$\Box C_{10}$. \label{fig:five}}

\begin{rmk}
The corollary above and the cases computed in Section~\ref{sec:small} and in~\cite{Pavlic2013} give that $\gamma_R(P_m\Box C_{n}) =\frac{2(m+1)n}{5}$ if $m,n\geq 4$ and $n\equiv 0\pmod 5$. Meanwhile, the Roman domination number is smaller if $m=2,3$.
\end{rmk}

\section{Conclusions}
{Algorithms}~\ref{alg:small} and~\ref{alg:recurrence} provide the exact value of $\gamma _R(P_m\Box C_n)$ where $m$ is fixed and their high computational requirements  restrict the maximum value of $m$ for which the algorithms are useful. In this work we have computed $\gamma _R(P_m\Box C_n)$ for $m=7,8,9$ according to the available platforms and these results could be extended by using more powerful platforms.
Additionally, we have computed the Roman domination number of cylinders $P_m\Box C_n$ in cases $m\geq$ {$10$}, $n\equiv 0\pmod 5$, expanding the family of such graphs whose Roman domination number is known. We have also provided a lower bound of this parameter if $m\geq$ {$10$} and $n\not\equiv 0\pmod 5$. We have implemented algorithms that allow to compute the exact values for selected cases and the general lower bound, by using the $(\min ,+)$ powers of large matrices in GPU platforms.

Our technique involves the computation of the loss of Roman dominating sets. Computing the loss has provided the exact values of several domination parameters in the Cartesian product of two paths (see~\cite{Goncalves2011,Rao2019}). We have applied similar techniques to the Roman domination number in cylinders and we have obtained the exact values only in some cases, but a lower bound in other ones. This different behaviour is due to several reasons. On the one hand, this technique computes the loss located at the border of the graph. Therefore, in the cylinder $P_m\Box C_n$ just upper and lower borders play a role, showing information about the loss that depends on the parameter $n$ but avoiding what depends on $m$. However, in $P_m\Box P_n$ the loss depending on both $m$ and $n$ can be computed by using the four borders.

On the other hand, we think that the formula for the general case will depend on the parity of $n$ $\pmod 5$. Meanwhile, the Cartesian product of two paths presents a unique formula for big enough $m$ and $n$. It is likely that a technique other than the computation of the loss will be necessary to completely solve this problem.

Moreover, it would be interesting to compute the exact values in additional small cases, in order to conjecture formulas for any value of both $m$ and $n$.

\appendices
\section{Proof of Lemma 10}
\label{appendix}

Our strategy to prove this lemma is similar to~\cite{Guichar2004}, so we are going to compute the loss just in both borders of the cylinder.

Denote $V(P_m)=\{u_1, u_2,\dots u_m\}$, $V(C_n)=\{v_1, v_2,\dots v_m\}$ and consider the following partition of $V(P_m\Box C_n)$:
$V_1=\{(u_i,v_j)\colon 1\leq i\leq 4 \text{ and } 1\leq j \leq n\}$, $V_2=\{(u_i,v_j)\colon 5\leq i\leq m-4 \text{ and }  1\leq j\leq n\}$ and $V_3=\{(u_i,v_j)\colon m-3\leq i\leq m \text{ and }  1\leq j\leq n\}$ (see Figure~\ref{fig:partition}). We call $G_i$ the subgraph of $P_m\Box C_n$ induced by the vertex set $V_i$, so $G_1$ and $G_3$ are both isomorphic to the cylinder $P_4\Box C_n$. Meanwhile, $G_2$ is isomorphic to $P_{m-8}\Box C_n$. If $f$ is a Roman dominating function, we denote by $f_i=f|_{V_i}$, the restriction of $f$ to the set $V_i$. Note that $f(P_m\Box C_n)=f_1(G_1)+f_2(G_2)+f_3(G_3)$ (the functions $f_i$ are not necessarily Roman dominating functions, but the concept of weight can obviously be extended to them).

\Figure[!h]()[width=0.18\textwidth]{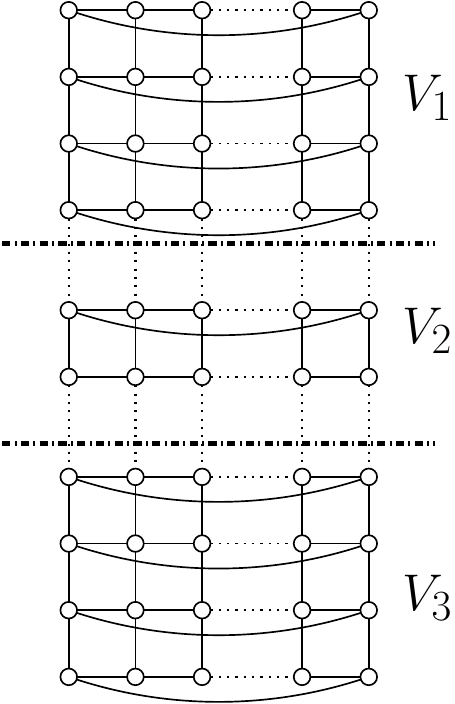}
{Partition of $V(P_m\Box C_n), m,n\geq${$10$}. \label{fig:partition}}

For $i=1,2,3$ we denote $D_i=(S_1\cap V_i)\cup N[(S_2\cap V_i)]$, which is the set of vertices dominated by $(S_1\cup S_2)\cap V_i$. Note that the closed neighborhood is $N[(S_2\cap V_i)]$ computed in the cylinder $P_m\Box C_n$, so it is not necessarily contained in $V_i$, and

\begin{equation}\label{eq:partition}
|D_i| \leq |(S_1\cap V_i)|+ 5|(S_2\cap V_i)|
\end{equation}

Moreover, it is clear that $D_1\cup D_2\cup D_3=V(P_m\Box C_n)$ and $D_1,D_2,D_3$ are not necessarily disjoint.

Using this notation, we can rewrite the loss of a Roman dominating function as follows:\\
$L(f)=\frac{5}{2}f(P_m\Box C_n)-mn=\frac{5}{2}(f_1(G_1)+f_2(G_2)+f_3(G_3))-|D_1\cup D_2\cup D_3|\geq (\frac{5}{2}f_1(G_1)-|D_1|)+(\frac{5}{2}f_2(G_2)-|D_2|)+(\frac{5}{2}f_3(G_3)-|D_3|)\geq (\frac{5}{2}f_1(G_1)-|D_1|)+(\frac{5}{2}f_3(G_3)-|D_3|)$.

The first inequality comes from the sets $D_i$ not being necessarily disjoint and the second one comes from the fact that each term of the sum is positive or zero, by using the inequality shown in Equation~\ref{eq:partition}:
$
\renewcommand*{\arraystretch}{1.2}
\begin{array}{l}
\frac{5}{2}f_i(G_i)-|D_i|=\frac{5}{2}(|S_1\cap V_i|+2|S_2\cap V_i|) -|D_i|\\
  \geq \frac{5}{2}(|S_1\cap V_i|+2|S_2\cap V_i|)-(|S_1\cap V_i|+ 5|S_2\cap V_i|)\\
  =\frac{3}{2}|S_1\cap V_i|.
\end{array}
$
In fact, we expect $|S_1\cap V_2|$ to be zero or close to zero in optimal cases, and this is why we discard the second term of the sum.

Vertex sets $V_1$ and $V_3$ play a similar role, so we will focus on $V_1$. The function $f_1$ satisfies that every vertex $v\in V_1$, not in row number four, with $f_1(v)=0$ is adjacent to at least a vertex $w$ satisfying $f_1(w)=2$. Meanwhile, vertices in the fourth row do not need to be dominated by vertices in $V_1$. This leads us to the following definition:

\begin{defn}
Consider $P_4\Box C_n$ as the subgraph consisting of the fourth top (or bottom) rows of $P_m\Box C_n, m\geq$ {$10$}. An almost Roman dominating function on $P_4\Box C_n$ is $g\colon V(P_4\Box C_n)\to \{0,1,2\}$  such that every vertex $v$, not in the fourth row, with $f(v)=0$ is adjacent to at least a vertex $w$ satisfying $f(w)=2$. We denote by $R^g_i=\{v\in V(P_4\Box C_n)\colon g(v)=i\}$ (we will omit $g$ if no confusion is possible).
\end{defn}

The set of dominated vertices of an almost Roman dominating function is $D(g)=R_1\cup N[R_2]$, where the closed neighborhood is computed in $P_m\Box C_n$. The loss of an almost dominating function is $L_a(g)=\frac{5}{2}g(P_4\Box C_n)- |D(g)|$ and the minimum loss of an almost Roman dominating function is $\mathcal{L}_a(n)=\min_g L_a(g)=\min_g (\frac{5}{2}g(P_4\Box C_n)- |D(g)|)$.

With this notation, the relationship between the loss of a Roman dominating function of $P_m\Box C_n$ and the loss in the borders of the cylinder is as follows:\\
$L(f)\geq  (\frac{5}{2}f_1(G_1)-|D_1|)+(\frac{5}{2}f_3(G_3)-|D_3|)\geq 2\mathcal{L}_a(n)$, therefore the minimum loss satisfies $\mathcal{L}(m,n)\geq 2\mathcal{L}_a(n)=\min_g (5g(P_4\Box C_n)- 2|D(g)|)$. We now modify  Algorithm~\ref{alg:small} and Algorithm~\ref{alg:recurrence}, in order to compute the value of $\min_g (5g(P_4\Box C_n)- 2|D(g)|)$.

We keep the definition of correct words, words of length $4$ in the alphabet $\{a,b,c,d\}$ not containing $ac, ca, ab, ba, bb$  and we need to modify the conditions for a word $\mathbf{p}=(p_1,p_2,p_3,p_4)$ to follow another word $\mathbf{q}=(q_1,q_2,q_3,q_4)$. First and intermediate row conditions are the same, but the fourth row, which is the last one in this case, behaves in a different way:

\begin{enumerate}
\item conditions for the first row: the same as before

\item conditions for the intermediate rows $2\leq i\leq 3$: the same as before

\item conditions for the last row (fourth row):
\begin{itemize}
\item if $q_4=a$, then $p_4=a$ or $p_4=c$
\item if $q_4=b$, then $\big(p_4=c \text{{ and} } p_{3}=a\big)$ or $p_3=d$
\item if $q_4=c$, then $p_4=a$ or $p_4=b$ or $\big(p_4=c \text{{ and} } p_{3}=a\big)$ or $p_4=d$
\item if $q_4=d$, then $p_4=a$ or $p_4=b$ or $\big( p_4=c \text{{ and} } p_{3}=a\big)$ or $p_4=d$
\end{itemize}
\end{enumerate}

Just conditions in case $q_4=d$ are new, because vertices in the fourth row do not need to be dominated. Therefore, in this case $p_4=a$ is possible but not compulsory and other labels for $p_4$ are suitable.

Every almost Roman dominating function is uniquely identified with an ordered list of correct words such that each word can follow the previous one and the first word can follow the last one.

The directed graph $G$ is defined as before: vertices are the correct words and there is an arc from $\mathbf{q}$ to $\mathbf{p}$ if $\mathbf{p}$ can follow $\mathbf{q}$. We need an appropriate labeling $\ell $ of the arcs such that, $\ell (Q)=5g(P_4\Box C_n)- 2|D(g)|$, for $Q=(\mathbf{p_1 p_2})\dots (\mathbf{p_{n-1} p_{n}})(\mathbf{p_n p_1})$ a closed path in $G$ and $g$ the almost Roman dominating function represented by $Q$.

For an arc $(\mathbf{qp})$, we define $\ell(\mathbf{q,p})=10\mathbf{p}(a)+5\mathbf{p}(b)-2nd(\mathbf{q,p})$, where $\mathbf{p}(a)$ is the number of $a's$ of $\mathbf{p}$, $\mathbf{p}(b)$ is the number of $b's$ of $\mathbf{p}$ and $nd(\mathbf{q,p})$ is the number of newly dominated vertices (see~\cite{Guichar2004}). That is, the number of vertices dominated by $\mathbf{p}$ which are not dominated by $\mathbf{q}$, when $\mathbf{q,p}$ are consecutive columns of an almost Roman dominating function. We compute such number with Algorithm~\ref{alg:newly_dominated}.

\begin{algorithm}[H]
\caption{Computation of $nd(\mathbf{q,p})$}
\label{alg:newly_dominated}
\begin{algorithmic}[1]
\REQUIRE {$\mathbf{q,p}$ correct words such that $\mathbf{p}$ can follow $\mathbf{q}$}
\ENSURE {$nd(\mathbf{q,p})$}
\STATE $nd(\mathbf{q,p})=0$
\FOR {$i=1$ to $4$}
\SWITCH {$q_i, p_i$}
\CASE{$q_i=a \text{ and } p_i=a$\,}
 \STATE $nd(\mathbf{q,p})=nd(\mathbf{q,p})+1$
\ENDCASE
\CASE{$q_i=b \text{ and }  p_i=c$\,}
 \STATE $nd(\mathbf{q,p})=nd(\mathbf{q,p})+1$
\ENDCASE
\CASE{$q_i=c \text{ and }  p_i=a$\,}
 \STATE $nd(\mathbf{q,p})=nd(\mathbf{q,p})+2$
\ENDCASE
\CASE{$q_i=c,d \text{ and } p_i=b,c$\,}
 \STATE $nd(\mathbf{q,p})=nd(\mathbf{q,p})+1$
\ENDCASE
\CASE{$q_i=d \text{ and }  p_i=a$\,}
 \STATE $nd(\mathbf{q,p})=nd(\mathbf{q,p})+3$
\ENDCASE
\ENDSWITCH
\ENDFOR
\IF {$p_4=a$}
\STATE $nd(\mathbf{q,p})=nd(\mathbf{q,p})+1$
\ENDIF
\RETURN $nd(\mathbf{q,p})$
\end{algorithmic}
\end{algorithm}

Note that the computational complexity of Algorithm~\ref{alg:newly_dominated} is of the order of the square of the number of correct words, so it does not increase the complexity of the computation of the matrix $A(G)$.

If $g$ is the almost dominating function represented by the closed path $Q=(\mathbf{p_1 p_2})(\mathbf{p_2 p_3})\dots (\mathbf{p_{n-1} p_n})
(\mathbf{p_n p_1})$ then $\Big(\sum _{i=1}^{n-1}nd (\mathbf{p_i, p_{i+1}})\Big) +nd(\mathbf{p_n p_1})=|D(g)|$, the number of vertices dominated by $g$. Therefore, we obtain that $\ell (Q)=\Big(\sum _{i=1}^{n-1}\ell(\mathbf{p_{i} p_{i+1}})\Big) +\ell(\mathbf{p_n p_1})=\Big(\sum _{i=1}^{n-1} 10\mathbf{p_{i+1}}(a)+5\mathbf{p_{i+1}}(b)-2nd(\mathbf{p_i,p_{i+1}})\Big)+10\mathbf{p_1}(a)+5\mathbf{p_1}(b)-2nd(\mathbf{p_n,p_1})= 10\sum _{i=1}^{n} \mathbf{p_i}(a)+ 5\sum _{i=1}^{n} \mathbf{p_i}(b)-2\Big(\big(\sum _{i=1}^{n-1}nd (\mathbf{p_i, p_{i+1}})\big) +
nd(\mathbf{p_n p_1})\Big)=5g(P_4\Box C_n)- 2|D(g)|$,

The construction of matrix $A(G)$ is as follows:
$$A(G)_{\mathbf{qp}} =
\left\{
\begin{array}{ll}
\ell(\mathbf{q,p}) & \text{if }(\mathbf{q p}) \text{ is an arc of } G,\\
\infty & \text{otherwise.}
\end{array}
\right.
$$

Algorithm~\ref{alg:recurrence} for $m=4$, using the new rules for construction of the matrix $A(G)$, gives $n_0=30$, $a=1$, $b=1$. This means $2\mathcal{L}_a(n) =\min_g (5g(P_4\Box C_n)- 2|D(g)|)$ satisfies the finite difference equation $ 2\mathcal{L}_a(n+1)- 2\mathcal{L}_a(n)=1$, for $n\geq 30$. Moreover, by using Algorithm~\ref{alg:small} (again with the new matrix $A(G)$), we have obtained $2\mathcal{L}_a(n) =n$, for {$10$}$\leq n\leq 30$. This gives $2\mathcal{L}_a(n)=n$, for $n\geq${$10$}, and finally, $\mathcal{L}(m,n)\geq 2\mathcal{L}_a(n)=n$, for $m,n\geq${$10$}, as desired.

\EOD
\end{document}